\documentclass[12pt]{article}
\usepackage{graphics}
\usepackage{epsfig}
\usepackage{amsthm, amsmath,amssymb,amsfonts}

\usepackage{color}

\usepackage{hyperref}
\hypersetup{colorlinks, citecolor=blue, filecolor=blue, linkcolor=blue, urlcolor=blue}

\newtheorem{theorem}{Theorem}[section]
\newtheorem{proposition}[theorem]{Proposition}
\newtheorem{corollary}[theorem]{Corollary}
\newtheorem{lemma}[theorem]{Lemma}
\newtheorem{remark}[theorem]{Remark}
\newtheorem{definition}[theorem]{Definition}
\newtheorem{example}[theorem]{Example}

\def\res{{\rm Res}}
\def\pp{{\rm pp}}
\def\gg{\mathfrak{g}}

\def\cc{{\cal C}}
\def\ee{{\cal E}}
\def\adj{{\cal A}_{d-2}(\ccc)}

\def\jadj{{\cal A}_{d-3}(\ccc)}

\def\genuss{{\rm\bf genus}}

\def\genus{{\rm genus}}
\def\sing{{\rm sing}}
\def\ccc{\overline{\cal C}}

\def\K{\mathbb{K}}
\def\F{\mathbb{F}}
\def\E{\mathbb{E}}
\def\L{\mathbb{L}}
\def\pj{\mathbb{P}}

\def\cP{{\cal P}}

\begin{document}

\title{Radical Parametrization of Algebraic Curves by Adjoint Curves}
\author{
J. Rafael Sendra\footnote{Supported by the Spanish ``Ministerio de
Educaci\'on e Innovaci\'on" under the Project MTM2008-04699-C03-01} \\
Dep. de Matem\'aticas \\
Universidad de Alcal\'a \\
Alcal\'a de Henares, Madrid, Spain \\
 {\tt rafael.sendra@uah.es}
\and
David Sevilla\thanks{Supported by the Spanish ``Ministerio de Educaci\'on e Innovaci\'on" under the Project MTM2007-67088} \\
Johann Radon Institute (RICAM) \\
Altenbergerstrasse 69 \\
A-4040 Linz, Austria \\
\texttt{david.sevilla@oeaw.ac.at}
}
\date{}
\maketitle
\begin{abstract}
We present algorithms for parametrizing by radicals an irreducible curve, not necessarily plane, when the genus is less o equal to 4 and they are defined over an algebraically
closed field of characteristic zero. In addition, we also present an algorithm for parametrizing by radicals any irreducible plane curve of degree $d$
having at least a point of multiplicity $d-r$, with $1\leq r \leq 4$ and, as a consequence, every irreducible plane curve of degree $d \leq 5$ and every irreducible singular plane curve of degree 6.

\end{abstract}



\section{Introduction}

It is well known that the only algebraic curves parametrizable by rational parametrizations are those of genus zero, and there are algorithms for that purpose (see e.g. \cite{PerezSendraWinkler2008a}). However, in many applications, this is an strong limitation because either the curves appearing in the process are not rational (i.e. genus zero curves) or the algebraic manipulation of the geometric object does not preserve the genus; this happens, for instance, when applying offsetting constructions (see \cite{ArrondoSendraSendra1997a}) or performing conchoidal transformations (see \cite{SendraSendra2010a}).

Motivated by this fact, we analyze in this paper the problem of developing algorithms to parametrize (not rationally) a bigger class of algebraic curves. In this sense, we consider the radical parametrizations. Essentially, a radical parametrization is a family of finitely many parametrizations given by rational functions whose numerators and denominators are radicals expression of polynomials (see Def. \ref{def-radical-parametrization} for a formal definition). For instance, the curve given by the polynomial $x^4+y^4-1$ cannot be parametrized rationally but admits the radical parametrization
\[ \{ (\pm\sqrt[4]{-1-t^4}, t) \} \cup \{ (\pm\sqrt{-1} \sqrt[4]{-1-t^4}, t) \}. \]

In \cite{Zariski1926a}, Zariski proved that the general complex projective curve of genus $\gg>6$ is not parametrizable by radicals. Moreover, as remarked in \cite{PirolaSchlesinger2005a}, Zariski's result is sharp. Indeed, a result within the Brill-Noether theory (see \cite{BrillNoether1873a}, or \cite[Chapter V]{ArbarelloCornalbaGriffithsHarris1985a} for a more modern account) states that a curve of genus $\gg$ has a linear system of dimension $1$ and degree $\lceil\frac{\gg}{2}+1\rceil$ \cite[p. 206]{ArbarelloCornalbaGriffithsHarris1985a}, thus a map of that degree to $\pj^1$ (and possibly lower degree maps as well). Therefore for $\gg=3,4$ there exists a $3:1$ map whose inversion would provide a radical parametrization with cubic roots, and for $\gg=5,6$ the inversion of the existing $4:1$ would provide a radical parametrization with quartic roots.


In this paper we present algorithms for computing radical parametrizations of irreducible, not necessarily plane, algebraic curves when the genus is less or equal to $4$ and they are defined over an algebraically closed field of characteristic zero. Also, we prove that although offsets of rational curves are not necessarily rational, offsets
of curves parametrizable by radicals are parametrizable by radicals (see Prop. \ref{prop-offsets}). A similar result is stated for conchoids (see Prop. \ref{prop-conchoids}). 

In Section 2 we introduce the notion of radical parametrization and we state some preliminary results. Next, in Section \ref{sec-lines}, we see that the classical idea of parametrizing rationally using lines can be extended to this new context, provided that there exists a point of multiplicity $d-r$, where $d$ is the degree of the curve and $1\leq r \leq 4$; note that $r=1$ corresponds to monomial curves. As a consequence, every irreducible plane curve of degree less or equal $5$ (take any point on the curve), and every singular irreducible curve of degree 6 (take a singular point) are parametrizable by radicals. Finally, in Section \ref{sec-adjoints}, we provide algorithms, based on linear systems of adjoints curves, for parametrizing by radicals every irreducible curve (not necessarily plane) of genus less or equal 4.

We finish this introduction remarking that, although we do not provide algorithms for parametrizing by radicals genus 5 and 6 curves, we provide algorithms for the cases of genus $\leq$ 4 and for $d$ degree curves having a point of multiplicity $d-r$, where $1\leq r \leq 4$. This implies a clear increase of the family of curves for which we can compute a parametrization, either rational or radical. We leave for future research the cases of genus 5 and 6.

\section{Radical Parametrizations}\label{sec-radicals}

We start by recalling some preliminary concepts from Galois theory; for further details see e.g. \cite{Cohn2003a}. Let $K$ be a field of characteristic zero. We say that $f\in K[x]$ is {\sf soluble or solvable by radicals over $K$} if there exists a finite tower of field extensions
\[ K=K_0\subset K_1\subset \cdots \subset K_r\]
such that
\begin{enumerate}
\item for $i=1,\ldots,r$, $K_i/K_{i-1}$ is the splitting field of a polynomial of the form $g_i(x)=x^{\ell_i}-a_i\in K_{i-1}[x]$, for $\ell_i>0$ and $a_i\neq 0$;
\item the splitting field of $f$ over $K$ is contained in $K_r$.
\end{enumerate}
A tower as above is called a {\sf root tower for $f$ over $K$}.

A central theorem in the theory states that $f\in K[x]$ is solvable by radicals iff its Galois group is solvable. Now,
let $f(t)$ be the general equation of degree $n$ over $K$, i.e. $f(x)=x^n-y_1x^{n-1}+\cdots+(-1)^n y_n\in K(y_1,\ldots,y_n)[x]$, where $y_i$ are unknowns.
Then, the theorem of Abel states that the Galois group of $f(x)$ is the symmetric group $S_n$, and hence $f(x)$ is solvable by radicals iff $n\leq 4$.

Now, and throughout the rest of the paper, let $\F$ be an algebraically closed field of characteristic zero,
$t$ a transcendental element over $\F$ and $\K=\F(t)$. In the following, we introduce the notions of square parametrization and radical parametrization. Essentially, an square parametrization is a parametrization given by rational functions which numerators and denominators are radicals expression of polynomials, and a radical parametrization is a finite set of square parametrizations. More precisely, we have the next definition.



\begin{definition}\label{def-square-parametrization}
 $(R_1,\ldots,R_n)$ is called an {\sf (affine) square parametrization} if there exists $g(x)\in \K[x]$,
soluble over $\K$, such that $(R_1,\ldots,R_n)\in \E^{n}\setminus \F^n$, where $\E$ is the last field extension of
a root tower of $g$ over $\K$.

Similarly,
 $(R_1:\cdots:R_n:R_{n+1})\in \pj^n(\E) \setminus \pj^n(\F)$ is called a {\sf (projective) square parametrization}.
\end{definition}

\begin{lemma}\label{lema-def-square-param}
Let $\K=\F(t)=\K_0\subset \K_1\subset \cdots \subset \K_r=\E$ be a root tower for $g(x)\in \K[x]$ over $\K$, and $S\in \E$.
Then, for all $t_0\in \F$, but finitely many exceptions, $S(t_0)$ is well defined.
\end{lemma}

\begin{proof} Let $g_{i}(x)=x^{\ell_i}-a_i\in \K_{i-1}[x]$, where $\ell_i>0$ and $a_i\neq 0$, be such that
$\K_i$ is the splitting field of $g_i$ over $\K_{i-1}$. We reason by induction over the field extension in the tower. If $S\in \K_0$ the claim in trivial.
 Let us assume that the statement is true over $\K_{i-1}$ and let $S\in \K_i$. Let $\alpha_i$ be a root of $g_i$.
Since $\F$ is algebraically closed then $\K_{i}=\K_{i-1}(\alpha_i)$ (see e.g. \cite[Prop. 7.10.7]{Cohn2003a}).
Then $S= \beta_0 +\beta_1 \alpha_i+\cdots +\beta_{\ell_i-1} \alpha_i^{\ell_i-1}$, with $\beta_j\in \K_{i-1}$. So, by induction,
 for almost all $t_0\in \F$, $\beta_j(t_0)$ is well defined; also $a_i(t_0)$, and since $\alpha_i^{\ell_i}=a_i$ and $\F$ is algebraically closed, $\alpha_{i}(t_0)$ too. Hence, also $S$.
\end{proof}

\begin{corollary}\label{cor-def-square-param}
Let $R$ be an affine parametrization. For all $t_0\in \F$, but finitely many exceptions, $R(t_0)$ is well defined.
\end{corollary}

\begin{definition}\label{def-radical-parametrization}
We say that an irreducible affine algebraic curve ${\cal E}\subset \F^n$ is {\sf soluble by radicals} or {\sf pa\-ra\-me\-tri\-za\-ble by radicals} if there exists a
finite set of square parametrizations $\{ \cP_i(t) \}_{i=1,\ldots,r}$ such that
\begin{enumerate}
\item for all $i\in \{1,\ldots,r\}$ and for all $t_0\in \F$ such that $\cP_i(t_0)$ is defined then $\cP_i(t_0)\in {\cal E}$;
\item for all but finitely many $P\in \cal E$ there exists $t_0\in \F$ and $i\in \{1,\ldots,r\}$ such that $P=\cP_i(t_0)$.
\end{enumerate}
In that case, $\{ \cP_i(t) \}_{i=1,\ldots,r}$ is called a {\sf radical (affine) parametrization of $\cal E$.}
\end{definition}

\begin{remark} We observe that
\begin{enumerate}
\item every rational parametrization can be seen as an square parametrization, and hence every rational curve is soluble by radicals,
\item considering projective square parametrizations, one introduces similarly the notion of being projectively parametrizable by radicals.
\end{enumerate}
\end{remark}


\begin{example} The Fermat cubic $x^3+y^3=1$ (which has genus 1) can be parametrized by radicals as
\[ \cP(t)=\{(\alpha_1\sqrt[3]{1-t^3},t),\ (\alpha_2 \sqrt[3]{1-t^3},t),\ (\alpha_3 \sqrt[3]{1-t^3},t) \}, \]
where $\alpha_1=1, \alpha_2=\frac{-1+\sqrt{3}\,i}{2}, \alpha_3=\frac{-1-\sqrt{3}\,i}{2}$; i.e. the roots of $x^3-1$.
First we see that each pair is a square parametrization by considering the extension sequence
\[ \mathbb{C}(t)=\K_0\subset \K_1=\K_0(\sqrt[3]{1-t^3})=\E \]
with $g_1(x)=x^3-(1+t^3)$. Moreover, (1), (2) in Definition \ref{def-radical-parametrization} clearly hold.
\end{example}

From Def. \ref{def-radical-parametrization} one deduces that the property of being soluble by radicals is preserved under birational transformations.

\begin{proposition}\label{prop-birational}
Let $\ee$ and $\ee^*$ be birationally equivalent curves over $\F$. Then $\ee$ is soluble by radicals iff $\ee^*$ is soluble by radicals.
\end{proposition}

Taking into account that every curve is birationally equivalent to a plane curve (see e.g. \cite[p. 155]{Fulton1969a}), we may work without loss of generality with plane algebraic curves. So, we introduce the following additional notation that will be used throughout the paper. Let $\cc$ be a plane irreducible affine curve over
$\F$ and $f(x_1,x_2)$ its defining polynomial. We denote by $\ccc$ the projective closure of $\cc$ and by $F(x_1,x_2,x_3)$ the homogenization
of $f(x_1,x_2)$. Moreover, we denote by $\deg(\ccc)$ or by $\deg(\cc)$ the degree of $\ccc$ (i.e. the degree of $F$),
by $\deg_{x_i}(\ccc)$ the partial degree of $F$ with respect to $x_i$, by $\genus(\cc)$ or by $\genus(\ccc)$ the genus of $\cc$,
and by $\sing(\ccc)$ the singular locus of $\ccc$.

As a first immediate result, the fact that a polynomial of degree at most $4$ in $\K[x]$ is soluble by radicals yields the following proposition.

\begin{proposition}\label{prop-grados-parciales}
If there exists $i\in \{1,2,3\}$ such that $\deg_{x_i}(\ccc)\leq 4$ (in particular if $\deg(\ccc)\leq 4$) then $\ccc$ is parametrizable by radicals.
\end{proposition}

\begin{proof} If $\cc$ is a line the result is trivial. Assume without loss of generality that $\deg_{x_1}(\ccc)=r\leq 4$. Define $g(x_1)=f(x_1,t)\in \K[x_1]$ which is
soluble over $\K$. Let $\E$ be the last field of a root tower for $g$ over $\K$, and $\{\alpha_1,\ldots,\alpha_r\}$ the roots of $g$ over $\K$;
we write $g$ as
\[ g(x_1)=\lambda(t) \prod_{i=1}^{r} (x_1-\alpha_i(t)). \] Then, for $i=1,\ldots,r$, $(\alpha_i(t),t)\in \E^2\setminus \F^2$.
So they are (affine) square parametrizations. Since $f(\alpha_i(t),t)=0$,
condition (1) in Def. \ref{def-radical-parametrization} holds. Now, let $\Delta\subset \F$
be a finite subset such that,
for $i=1,\ldots,r$ and for all $t_0\in \F\setminus \Delta$, $\alpha_i(t_0)$
is well defined (see Lemma \ref{cor-def-square-param}) and $\lambda(t_0)$ is well
defined and non-zero. Let $(a,b)\in \cc$ such that $b\not\in \Delta$. As $\cc$ is not a line and it is irreducible, we are excluding finitely many points on $\cc$.
So, since $\lambda(b)\neq 0$ then $\prod_{i}^{r} (a-\alpha_i(b))=0$. Hence, there exists $i$ such that $(\alpha_i(b),b)=(a,b)$.
\end{proof}

\begin{example}\label{example-grado-parcial}
We consider the curve $\cc$ of degree 14 defined by $$f(x_1,x_2)=x_{1}^4x_{2}^{10}+x_{2}x_{1}+1$$ whose genus is 1.
Since $\deg_{x_1}(\ccc)=4$ we apply Prop. \ref{prop-grados-parciales}, and solving by radicals the polynomial $f(x_1,t)\in \mathbb{C}(t)[x_1]$, we get the radical parametrization of $\cc$:
\[ \left\{\left({\frac {\sqrt {6}\sqrt {\Delta_{{3}}}\pm \sqrt {\Delta_{{4}}}}{{12 t
}^{3}}},t\right),\left(-\,{\frac {\sqrt {6}\sqrt {\Delta_{{3}}}\mp
\sqrt {\Delta_{{5}}}}{{12t}^{3}}},t\right)\right\}, \]
where
\[ \begin{array}{l}\Delta_1=-768\,{t}^{6}+81, \quad \Delta_2= 108+12\,\sqrt {\Delta_{{1}}}, \quad \displaystyle{\Delta_3={\frac {\sqrt[3]{\Delta_{{2}}^2}+48\,{t}^{2}}{\sqrt [3]{\Delta_{{2}}}}},} \\
\displaystyle{\Delta_4=-{\frac {6\,\sqrt {\Delta
_{{3}}}\sqrt[3]{\Delta_{{2}}^2}+288\,\sqrt {\Delta_{{3}}}{t}^{2}+72\,
\sqrt {6}\sqrt [3]{\Delta_{{2}}}}{\sqrt [3]{\Delta_{{2}}}\sqrt {\Delta
_{{3}}}}},} \\[1em]
\displaystyle{\Delta_5=-{\frac {6\,\sqrt {\Delta_{{3}}}\sqrt[3]{\Delta_{{2}}^2}+288\,
\sqrt {\Delta_{{3}}}{t}^{2}-72\,\sqrt {6}\sqrt [3]{\Delta_{{2}}}}{
\sqrt [3]{\Delta_{{2}}}\sqrt {\Delta_{{3}}}}}}.
\end{array}\]
\end{example}

\subsection{The Case of Offsets and Conchoids} 

In the introduction, we have mentioned that working only with genus zero curves can be a limitation in some applications are those requiring offsetting and conchoidal constructions. In this subsection, we see that both geometrical constructions are closed under radical parametrizations. We briefly recall here the intuitive idea of offset and conchoid; for further details we refer to \cite{ArrondoSendraSendra1997a} for the case of curves and \cite{SendraSendra2008a} for the case of conchoids. Let $\cc$ be an irreducible plane curve define by the polynomial $f(x,y)$, and $A=(\lambda_1,\lambda_2)$ a point in the plane, then (we recall that $\|(x,y)\|_{2}^{2}=x^2+y^2$).
\begin{enumerate}
\item the offset to $\cc$ at distance $d\in \F\setminus \{0\}$ is the Zariski closure of
\[ \{ P\pm \frac{d}{\|\nabla(f)(P)\|_2} \nabla(f)(P)\,\,\, \mbox{with $P\in \cc$ such that $|\nabla(f)(P)\|_2\neq 0$} \}\]
where $\nabla(f)(P)=(\frac{\partial f}{\partial x}(P),\frac{\partial f}{\partial y}(P))$.
\item the conchoid to $\cc$ at distance $d\in \F\setminus \{0\}$ from the focus $A$ is the Zariski closure of
\[ \{ P\pm \frac{d}{\|A-P\|_2} (A-P)\,\,\, \mbox{with $P\in \cc$ such that $|A-P\|_2\neq 0$} \}\]
\end{enumerate}
\begin{proposition}\label{prop-offsets}
Let $\cc$ be an irreducible plane curve parametrizable by radicals, and $d\in \F \setminus \{0\}$. The offset to $\cc$ at distance $d$ is also
parametrizable by radicals.
\end{proposition}

\begin{proof}
Let $\cP=\{{\cal P}_{i}(t)\}_{i=1,\ldots,r}$ be a radical (affine) parametrization of $\cc$. We consider the formal derivation with respect to $t$, and for each $\cP_i(t)=(R_1,R_2)\in \cP$ we define
$$O_{i}^{\pm}=\left(R_1\mp d \frac{(R_2)'}{\sqrt{(R_{1})'^{2}+(R_{2})'^{2}}},R_2\pm d \frac{(R_1)'}{\sqrt{(R_{1})'^{2}+(R_{2})'^{2}}}\right)$$
and $O=\{O_i^{\pm}(t)\}_{i=1,\ldots,r}$. Because of the definition of offset it is clear that $O$ and the offset to $\cc$ at distance $d$ satisfy conditions (1) and (2) in Def. \ref{def-radical-parametrization}. So, it only remains to proof that each $O_{i}^{\pm}$ is an square parametrization. To prove that, let $\cP_i=(R_1,R_2)$, and let
$\K_0=\F(t)\subset \cdots \subset \K_r$ be such that $R_1,R_2\in \K_r$. Reasoning by induction on $i$, one has that if $R\in \K_r$ then $R'\in \K_r$. Now, let $a=(R_{1})'^{2}+(R_{2})'^{2}\in \K_r$ and take $\K_{r+1}=\K_r(\sqrt{a})$. Then, $O_i^{\pm}\in \K_{r+1}^{2}\setminus \K_{0}^{2}$ and hence $O_i^{\pm}$ is an square parametrization.
\end{proof}

Reasoning similarly, one gets the following result.

\begin{proposition}\label{prop-conchoids}
Let $\cc$ be an irreducible plane curve parametrizable by radicals, $A\in \F^2$, and $d\in \F \setminus \{0\}$. The conchoid to $\cc$ from $A$ at distance $d$ is also
parametrizable by radicals.
\end{proposition}

\section{Radical parametrization by lines}\label{sec-lines}

Let us analyze more deeply the meaning of Prop. \ref{prop-grados-parciales}.
The fact that $\deg(\cc)\leq 4$ is intrinsic to the curve, while the property on the partial degrees is not.
In fact, most linear changes of coordinates modify the partial degrees while the fact of being soluble by radicals
is preserved (see Prop. \ref{prop-birational}).
So, we need a more geometrical interpretation of this fact.

Indeed, $\deg_{x_1}(\ccc)\leq 4$
iff $(1:0:0)$ is a $(d-r)$-fold point of $\ccc$, with $1\leq r\leq 4$; similarly with $(0:1:0)$, $(1:0:0)$ and $x_2, x_3$.
For instance, in Example \ref{example-grado-parcial}, $(1:0:0)$ has multiplicity $10$. So, if $\ccc$ has a singular point with this
multiplicity, performing a suitable change of coordinates we fulfill the hypothesis of Prop. \ref{prop-grados-parciales}.
This is equivalent to trying to parametrize (with radicals) using lines; note that if $r=1$, i.e. $\ccc$ is monomial,
we reach the case of rational curves parametrizable for lines (see \cite{PerezSendraWinkler2008a}, section 4.6).

So, let us assume without loss of generality that the origin is a $(d-r)$-fold point of $\cc$, with $1\leq r\leq 4$. Then,
\[ f(x_1,x_2)=f_{d}(x_1,x_2)+\cdots +f_{d-r}(x_1,x_2), \]
where $f_i$ is a homogeneous polynomial of degree $i$. Let $h(x_1,x_2)=x_1-t x_2$ be the defining polynomial of a pencil of
lines ${\cal H}(t)$ passing through the origin. Now, by B\'ezout's Theorem, ${\cal H}(t)$ intersects $\cc$ at the origin, with multiplicity $d-r$, and at $r$ additional points that depend on the slope $t$:
\[ f(x_1, t x_1)= x_{1}^{d-r} ( x_{1}^{r} f_d(1,t)+\cdots +f_{d-r}(1,t)). \]
Let $g(x_1)= f_d(1,t)x_{1}^{r}+\cdots +f_{d-r}(1,t)\in \K[x_1]$. By hypothesis $\deg(g)\leq 4$. So, solving $g(x_1)$ by radicals over $\K$
(say that $\alpha_i(t)$, with $i=1,\ldots,r$ are the roots), we get that
\[ \{ (\alpha_i(t),t\alpha_i(t)) \}_{i=1,\ldots,r}, \]
 is a radical parametrization of $\cc$.

This reasoning shows how to algorithmically parametrize by radicals the following family of curves.

\begin{theorem}\label{teorema-por-rectas}
If $\ccc$ has an $(d-r)$-fold point, with $d=\deg(\ccc)$ and $1\leq r\leq 4$, then $\ccc$ can be parametrized by radicals using lines.
\end{theorem}
\begin{proof}
Make a linear change of coordinates that moves the singularity to $(1:0:0)$ and apply Prop. \ref{prop-grados-parciales}.
\end{proof}

\begin{corollary}
Every irreducible plane curve of degree less or equal 5 is soluble by radicals.
\end{corollary}
\begin{proof}
Every point satisfies the condition of Theorem \ref{teorema-por-rectas}.
\end{proof}

\begin{corollary}
Every irreducible singular plane curve of degree less or equal 6 is soluble by radicals.
\end{corollary}
\begin{proof}
Take any singular point and apply Theorem \ref{teorema-por-rectas}.
\end{proof}

\begin{example}
Let $\cc$ be the curve, of degree $10$, defined by $$f(x_1,x_2)=x_{1}^{10}+x_{2}^{10}+x_{2}^{6}-2\,x_{1}^{6}.$$ The genus of $\cc$ is $21$ and the origin is a $6$-fold point. So, by Theorem \ref{teorema-por-rectas}, $\cc$ is parametrizable by radicals using lines. In this case $g(x_1)=x_{1}^4+t^{10}x_{1}^4+t^6-2$, which provides the radical parametrization
 \[ \left\{\left(
 \xi^i\sqrt[4]{\frac{-(t^6-2)}{1+t^{10}}},\xi^i t\sqrt[4]{\frac{-(t^6-2)}{1+t^{10}}}
\right)\right\}_{i=1,\ldots,4}, \qquad \xi^4=1.
\]

\end{example}

\section{Radical parametrization by adjoints curves}\label{sec-adjoints}

In this section, we see that curves of genus at most 4 can be parametrized by radicals. Obviously genus zero curves (i.e. rational) can be
parametrized by radicals. We show that following the adjoint curves approach for parametrizing rational curves
(see \cite{PerezSendraWinkler2008a}, section 4.7), a method for parametrizing by radicals curves of genus less or equal $3$ is derived.
In the second part we see how to extend the method for genus 4 curves using adjoints of lower degree.

For this purpose, in the sequel, $d=\deg(\ccc)$, $\gg=\genus(\ccc)$, and
${\cal A}_i(\ccc)$ is the linear system of adjoint curves to $\ccc$ of degree $i$. Also for an effective divisor $D$,
we denote by ${\cal H}(i,D)$ the linear system of curves of degree $i$ generated by $D$ (see sections 2.5 and 4.7. in \cite{PerezSendraWinkler2008a}).

\subsection{The case $\mathbf{\genuss(\ccc)\leq 3}$}

 In \cite{Noether1884a} it is proved that for $i\geq d-3$, the linear conditions for the $i$-degree adjoints, derived from the singularities of $\ccc$, are independent.
Now, the genus formula states (see \cite[Theorem 3.11]{PerezSendraWinkler2008a}) that
\[
 \gg = \frac{(d-1)(d-2)}{2} - \sum m_P(m_P-1).
\]
where sum applies to all points on the curve including the neighboring ones. From this we obtain that $\dim(\adj)=d-2+\gg$.

Now we choose $\ell$ different simple points
$\{Q_1,\ldots,Q_\ell\}$ in $\ccc$ such that $$\adj^*=\adj\cap {\cal H}\left(d-2,\sum_{i=1}^{\ell} Q_i\right)$$
has dimension $1$. Note that $\ell \geq (d-3)+\gg$. The defining polynomial of $\adj^*$ can be written as
$H^*(x_1,x_2,x_3,t)=\Phi_1(x_1,x_2,x_3)+t\,\Phi_2(x_1,x_2,x_3)$,
where $\Phi_1,\Phi_2$ are defining polynomials of adjoint curves in $\adj^*$.

Let us see $\ccc$ and $\adj^*$ as curves in $\pj^3(\K)$.
Since $\ccc$ is irreducible and $\adj^*$ has smaller degree, by B\'ezout's Theorem, $\ccc \cap \adj^*$ contains (counted with multiplicity)
$d(d-2)$ points. On the other hand, ${\cal B}=\sing(\ccc)\cup \{Q_1,\ldots,Q_{\ell}\}\subset \ccc \cap \adj^*$.
 Moreover,
because of the genus formula
, points in $\cal B$ count in that intersection
at least $s=(d-1)(d-2)-2\gg+\ell$. Therefore, at most
$d(d-2)-s$
intersections (counted with multiplicity) are free. In addition, since all points in $\cal B$ are defined over $\F$, at most these $d(d-2)-s$
intersection points are in $\pj^3(\K)\setminus \pj^3(\K)$. Finally, observe that
\[ d(d-2) - s=(d-2) + 2\gg - \ell\leq \gg+1. \]
In this situation, let $\res_w$ denote the resultant with respect to the variable $w$, and $\pp_w$ its primitive part with respect to $w$.
We consider the homogeneous polynomials
\[ S_1(x_1,x_3,t)=\pp_t(\res_{x_2}(F,H^*))\in \K[x_1,x_3], \] \[ S_2(x_2,x_3,t)=\pp_t(\res_{x_1}(F,H^*))\in \K[x_2,x_3]. \]
The linear factors of $S_1,S_2$ over $\K$ provide the points in $\ccc\cap \adj^* \cap \pj^2(\K)\setminus \pj^2(\F)$.
So, dehomogenizing with respect to $x_3$, say $R_i=S_i(x_i,1,t)$, we get that $R_1,R_2$ describe the $x_1$ and $x_2$ coordinates of the affine intersections of
$\cc\cap \adj^*$ in $\K^2\setminus \F^2$. Therefore, if $\delta_i=\deg_{x_i}(R_i)$ then
 $$1\leq \delta_i \leq (d-2)+2\gg-\ell\leq \gg+1.$$
Furthermore, for every root $\alpha_i$ of $R_1$ in $\K$, there exist $\{\beta_1,\ldots,\beta_{k(i)}\}$, roots of $R_2$, such that $f(\alpha_i,\beta_{k(i)})=0$.

Now, if $\gg\leq 3$ then $\delta_i\leq 4$. So, the roots $\{\alpha_i\}_{1\leq i \leq \delta_1}$ and
$\{\beta_j\}_{1\leq j \leq \delta_2}$ can be computed by radicals over $\K$. In the next theorem we prove that
\[ \cP=\{ \cP_{i,j}=(\alpha_i,\beta_{j})\}_{i=1,\ldots \delta_1,j=1\ldots,k(i)} \]
is, indeed, a radical parametrization of $\cc$. In order to actually compute $\cP_{i,j}$, i.e. the corresponding $\beta_j$ for each $\alpha_i$,
one may proceed as follows. One possibility is checking for which of the $\delta_1 \delta_2\leq 16$ substitutions one gets $f(\alpha_i,\beta_j)=0$.
A second possibility consists in computing, for each $\alpha_i$, a polynomial whose roots are the corresponding $\beta_j$. More precisely,
for each irreducible factor $T(x_1,t)$ over $\K$ of $R_1(x_1,t)$, let $\L_T$ be the quotient field of $\K[x_1]/(T)$, and
$$M_T(x_1,x_2,t)=\gcd_{\L_T[x_2]}(f(x_1,x_2),R_2(x_2,t)).$$
Note that $1\leq \deg_{x_i}(M_T)\leq \delta_i\leq 4$.
Then, for each root $\alpha_i$ of $T$ the corresponding $\beta_j$ are the roots in $\K$ of $M_T(\alpha_i,x_2,t)$.

This reasoning provides an algorithm to parametrize by radicals curves of genus $\leq 3$ using adjoint curves.

\begin{theorem}\label{theorem-genus-3}
Every irreducible curve $\cc$ (not necessarily plane) with $\genus(\cc)\leq 3$ can be parametrized by radicals using adjoint curves of degree $d-2$.
\end{theorem}
\begin{proof}
As we have already remarked, we may assume that the curve is plane.
We see that $\cP$ (see above) is an affine radical parametrization of $\cc$.
Since $\delta_i\leq 4$, $g(x)=R_1(x,t)R_2(x,t)\in \K[x]$ is soluble by radicals over $\K$. Let $\E$ be the last field in a
root tower for $g$ over $\K$. Then, $\cP_{i,j}\in \E^2\setminus \F^2$. So they are square parametrizations. In addition,
$f(\cP_{i,j})=0$. So condition (1) in Def. \ref{def-radical-parametrization} is satisfied. To prove condition (2)
in Def. \ref{def-radical-parametrization}, we proceed as follows: first we define a subset $\Sigma \subset \ccc$, secondly we prove that $\Sigma$ is
finite, and finally for each $P\in (\ccc \setminus \Sigma)$ we prove that there exists $t_0\in \F$ and $\cP_{i,j}$ such that $\cP_{i,j}(t_0)=P$.
For this purpose, for $\lambda\in \F$ we denote by $\mathfrak{I}(\lambda)$ the intersection of $\ccc$ and the adjoint $\cal D$ defined by
$\Phi_1+\lambda \Phi_2$. Note that since $\ccc$ is irreducible and $\deg(\ccc)> \deg({\cal D})$ then $\mathfrak{I}(\lambda)$ is finite.
In this situation, let $\Sigma$ be
the set of points $P\in \ccc$ such that at least one of the following holds:
\begin{enumerate}
\item[(i)] $P$ is not an affine point,
\item[(ii)] $\Phi_2(P)=0$,
\item[(iii)] $P\in \mathfrak{I}(\lambda)$, where $\lambda$ is such that either:
\begin{itemize}
\item for some $i=1,2$, the leading coefficient of $R_i$ with respect to $x_i$ is not defined at $\lambda$ or it is defined but vanishes at $\lambda$,
\item or $S_i(x_2,x_3,\lambda)$, or $S_2(F,H^*)(x_1,x_3,\lambda)$, or $\alpha_i(\lambda)$, for some $1\leq i\leq \delta_1$,
 or $\beta_j(\lambda)$, for some $1\leq j \leq \delta_2$, is not defined,
\end{itemize}
\item[(iv)] for $i=1,2$ the content with respect to $t$ of $\res_{x_i}(F,H^*)$ vanishes at $P$,
\item[(v)] For each irreducible factor $T$ of $R_1$ over $\K$ we consider the gcd polynomial $M_T$. Moreover, since $\L_T[x_2]$ is an Euclidean
domain, there exists
$U_T,V_T\in \L_T[x_2]$ such that $M_T=U_T f+V_T R_2$. Let $N_T(t)=\res_{x_1}(T,A)$ where $A$ is the leading coefficient of $M_T$ with respect to $x_2$; note
that since $A$ is not the zero of $\L_T$ then $N_T$ is not identically zero. \\
Then, $P\in \mathfrak{I}(\lambda)$, where $\lambda$ is such that either:
\begin{itemize}
\item $N_T$ is not defined at $\lambda$ or it is defined but $N_T(\lambda)=0$
\item or $U_T(x_1,x_2,\lambda)$ or $V_T(x_1,x_2,\lambda)$ or $M_T(x_1,x_2,\lambda)$ is not defined.
\end{itemize}
\end{enumerate}
Let us see that $\Sigma$ is finite. Since $\ccc$ is the projective closure of an affine curve, (i) provides finitely many points. Since
$\deg(\phi_2)=d-2$ and $F$ is irreducible of degree $d$, (ii) provides finitely many points. In (iii), because of the reasoning on
$\mathfrak{I}(\lambda)$ above, it suffices to show that we are considering finitely many $\lambda\in \F$. But this follows from Lemma \ref{lema-def-square-param},
and taking into account that the leading coefficient with respect to $x_i$
of $R_i$ and the coefficients of $S_i$ are in $\K=\F(t)$. So (iii) provides finitely many points.
About (iv), the content with respect to $t$ of $\res_{x_i}(F,H^*)$ is either
constant (in which case there is nothing to say) or a product of linear homogeneous polynomials. Since we may assume without loss of generality that $\ccc$
is not a line, the intersection points of these lines with $\ccc$ are finitely many. (v) follows as (iii).

Let $P=(a:b:1)\in (\ccc\setminus \Sigma)$ (see (i)). By (ii), $t_0=\Phi_1(P)/\Phi_2(P)$ is well defined. Since $F(P)=H^*(P,t_0)=0$, one has that
$\res_{x_2}(F,H^*)(a,1,t_0)=\res_{x_1}(F,H^*)(b,1,t_0)=0$ and, by (iv), $S_1(a,1,t_0)=S_2(b,1,t_0)=0$. Hence, $R_1(a,t_0)=R_2(b,t_0)=0$.
By (iii), the leading coefficient of $R_i$ with respect to $x_i$ does not vanish at $t_0$. So, again by (iii),
there exist $i$ and $j$ such that $a=\alpha_i(t_0)$ and $b=\beta_j(t_0)$. It remains to see that $j\in\{1,\ldots,k(i)\}$. Let
$T$ be the irreducible factor over $\K$ of $R_1$ that has $\alpha_i(t)$ as a root. Since $f(a,b)=R_2(b,t_0)=0$, and
$U_T(a,b,t_0),V_T(a,b,t_0),M_T(a,b,t_0)$ are well defined (see (v)), from the equality $M_T=U_Tf+V_T R_2$, one gets that $M_T(a,b,t_0)=0$. Moreover,
since $N(t_0)\neq 0$ (see (v)) and $T(a,t_0)=0$, then $A(a,t_0)\neq 0$. So, $j\in\{1,\ldots,k(i)\}$.
\end{proof}

\begin{example}
Let $\cc$ be the 10-degree irreducible curve over $\mathbb{C}$ defined by
\[\begin{array}{lcl} f(x_1,x_2) &= &{\frac {101977}{3375}}\,{x_{{1}}}^{5}+{\frac {11591}{125}}\,{x_{{1}}}^
{4}{x_{{2}}}^{2}-{\frac {263656}{3375}}\,{x_{{1}}}^{4}x_{{2}}+{x_{{2}}
}^{5}+{x_{{1}}}^{5}{x_{{2}}}^{5} \\ && \\ & & +{\frac {60781}{1125}}\,{x_{{1}}}^{3}{
x_{{2}}}^{2}-{\frac {515293}{6750}}\,{x_{{1}}}^{3}{x_{{2}}}^{3}+{x_{{1
}}}^{2}{x_{{2}}}^{4}+{\frac {164107}{6750}}\,{x_{{1}}}^{3}{x_{{2}}}^{4
}\\ && \\ && -{\frac {207329}{6750}}\,{x_{{1}}}^{5}x_{{2}}-{\frac {7967}{2250}}\,{
x_{{1}}}^{4}{x_{{2}}}^{5}-{\frac {88201}{6750}}\,{x_{{1}}}^{4}{x_{{2}}
}^{3}-\frac{5}{2}\,x_{{1}}{x_{{2}}}^{4}.\end{array} \]
 Its singularities are $(0:0:1), (1:0:0), (0:1:0)$, all three with multiplicity $5$, and
$(1:1:1)$ with multiplicity $3$. Note that Theorem \ref{teorema-por-rectas} is not applicable. However, $\genus(\ccc)=3$, so we can apply Theorem
\ref{theorem-genus-3}. To obtain $\adj^*$ we take $(d-3)+\gg=11$ simple points on $\ccc$. We consider $(-1:1:1),(1:2:1),(2:1:1), (-3:2:1), (-3:2:1),
(5:2:1),$ and $(-2\gamma:\gamma:1)$, where
\[ -{\frac {2973418}{1125}}+{\frac {1155472}{375}}\,\gamma-32\,{\gamma}^{
5}-{\frac {454012}{1125}}\,{\gamma}^{2}-{\frac {63736}{1125}}\,{\gamma
}^{4}=0. \]
The implicit equation of $\adj^*$ is

\vspace{1em}

 \noindent $H^*(x_1,x_2,x_3,t)=10\,{x_{{1}}}^{3}{x_{{2}}}^{2}{x_{{3}}}^{3}-2\,{x_{{1}}}^{4}{x_{{3}}}^
{4}-4\,{x_{{1}}}^{4}{x_{{2}}}^{2}{x_{{3}}}^{2}+{x_{{1}}}^{4}{x_{{2}}}^
{3}x_{{3}}-4\,{x_{{1}}}^{3}x_{{2}}{x_{{3}}}^{4}+5\,{x_{{1}}}^{4}x_{{2}
}{x_{{3}}}^{3}-8\,{x_{{1}}}^{3}{x_{{2}}}^{3}{x_{{3}}}^{2}+2\,{x_{{1}}}
^{3}{x_{{2}}}^{4}x_{{3}}+t (-3\,{x_{{1}}}^{2}{x_{{2}}}^{3}{x_{{3}}}^{3}+\frac{3}{2}\,{x_{{1}}}^{3}{x_{{2}}
}^{2}{x_{{3}}}^{3}-{x_{{1}}}^{4}{x_{{3}}}^{4}-\frac{1}{2}\,{x_{{1}}}^{4}{x_{{2
}}}^{2}{x_{{3}}}^{2}-{x_{{1}}}^{3}x_{{2}}{x_{{3}}}^{4}+\frac{3}{2}\,{x_{{1}}}^
{4}x_{{2}}{x_{{3}}}^{3}-\frac{1}{2}\,{x_{{1}}}^{3}{x_{{2}}}^{3}{x_{{3}}}^{2}+{
x_{{1}}}^{2}{x_{{2}}}^{4}{x_{{3}}}^{2}+2\,{x_{{1}}}^{2}{x_{{2}}}^{2}{x
_{{3}}}^{4}
)$

\vspace{1em}

\noindent and the polynomials $R_1$,$R_2$ are

\vspace{1em}

\noindent $R_1(x_1,t)=-6750\,{x_{{1}}}^{2}{t}^{5}-540000\,{x_{{1}}}^{2}{t}^{2}+2777264\,t{x_
{{1}}}^{2}-108000\,{x_{{1}}}^{2}-67500\,{x_{{1}}}^{2}{t}^{4}-270000\,{
x_{{1}}}^{2}{t}^{3}+4036120\,x_{{1}}{t}^{2}+36510\,x_{{1}}{t}^{4}+3651
\,x_{{1}}{t}^{5}+498844\,x_{{1}}{t}^{3}+108000\,x_{{1}}-162000\,x_{{1}
}t+336620\,{t}^{4}+270000\,t+216000+1323744\,{t}^{3}+31203\,{t}^{5},$

\vspace{1em}

\noindent $R_2(x_2,t)=47802\,{x_{{2}}}^{2}{t}^{4}+6750\,{x_{{2}}}^{2}{t}^{5}-216000\,{x_{{2}
}}^{2}+432000\,x_{{2}}+270000\,x_{{2}}t+656428\,x_{{2}}{t}^{3}-54000\,
x_{{2}}{t}^{2}+95604\,x_{{2}}{t}^{4}-3651\,x_{{2}}{t}^{5}-216000-31203
\,{t}^{5}-336620\,{t}^{4}-270000\,t-1323744\,{t}^{3}.
$

\vspace{1em}

Solving by radicals, and checking in this case where $f(\alpha_i,\beta_j)=0$, we get the radical parametrization
\[ \left\{ \left({\frac {108000+3651\,{t}^{5}+36510\,{t}^{4}+498844\,{t}^{3}+
4036120\,{t}^{2}-162000\,t\pm\sqrt {\Delta}t\pm2\,\sqrt {\Delta}}{2(135000\,{t}^{3}+33750\,{t}^{4}+3375\,{t}^{5}+270000\,{t}^{2}-
1388632\,t+54000)}},\right.\right. \]
\[ \left.\left. {\frac {-432000-270000\,t-656428\,{t}^{3}+
54000\,{t}^{2}-95604\,{t}^{4}+3651\,{t}^{5}\pm\sqrt {\Delta}t}{
12(7967\,{t}^{4}+1125\,{t}^{5}-36000)}}\right) \right\} \]
where

\vspace{1em}

\noindent $ \Delta=855810801\,{t}^{8}+14356902816\,{t}^{7}+104452412520\,{t}^{6}+
379019035776\,{t}^{5}+425890947184\,{t}^{4}+8076240000\,{t}^{3}+
190450224000\,{t}^{2}-605721024000\,t+26244000000.$
\end{example}

\subsection{The case $2\leq \mathbf{\genuss(\ccc)\leq 4}$}

Let $\gg\geq 2$, and let us consider now $\jadj$; i.e. adjoints of degree $d-3$. Then, $\dim(\jadj)=\gg-1$ similarly as above.
If $\gg>2$, we choose $\ell$ different simple points $\{Q_1,\ldots,Q_\ell\}$ in $\ccc$ such that $$\jadj^*=\jadj\cap {\cal H}\left(d-3,\sum_{i=1}^{\ell} Q_i\right)$$
has dimension $1$. Note that $\ell \geq \gg-2$. If $\gg=2$, then $\jadj=\jadj^*$. The defining polynomial of $\adj^*$ can be written as
$H^*(x_1,x_2,x_3,t)=\Phi_1(x_1,x_2,x_3)+t\,\Phi_2(x_1,x_2,x_3)$,
where $\Phi_1,\Phi_2$ are defining polynomials of adjoint curves in $\jadj^*$. Reasoning as above, we see $\ccc$ and $\jadj^*$
as curves in $\pj^3(\K)$. Then, $\sing(\ccc)\cup \{Q_1,\ldots,Q_{\ell}\}\subset \ccc \cap \jadj^*$.
 Therefore, at most $d(d-3)-(d-1)(d-2) +2\gg-\ell$
intersection points are in $\pj^3(\K)\setminus \pj^3(\F)$. But observe that
\[ d(d-3)-(d-1)(d-2) +2\gg-\ell =2(\gg-1)-\ell \leq 2(\gg-1)-\gg +2=\gg. \]
Therefore $R_1$ and $R_2$, defined as in the previous case, describe the $x_1$ and $x_2$ coordinates of the affine intersections of
$\cc\cap \jadj^*$ in $\K^2\setminus \F^2$. Moreover, if $\delta_i=\deg_{x_i}(R_i)$ then
 $$1\leq \delta_i \leq \gg.$$
Furthermore, for every root $\alpha_i$ of $R_1$ in $\K$, there exist $\{\beta_1,\ldots,\beta_{k(i)}\}$, roots of $R_2$,
such that $f(\alpha_i,\beta_{k(i)})=0$.
Now, if $2\leq \gg\leq 4$ then $\delta_i\leq 4$. Then reasoning as in the previous case we get the following theorem, which provides an algorithm to parametrize by radicals curves of genus $2,3$ and $4$ using adjoint curves.

\begin{theorem}\label{theorem-genus-4}
Every irreducible curve $\ee$ (not necessarily plane), with $2\leq \genus(\ee)\leq 4$, can be parametrized by radicals using adjoint curves of degree $(d-3)$.
\end{theorem}

\section{Conclusions}

The results presented in this paper provide algorithms for radical parametrizations of curves with a singularity of a high multiplicity, and also for curves of genus up to 4. Theoretical results mentioned in the introduction indicate that curves of genus 5 and 6 can also be parametrized by radicals. Considering adjoints of higher degree does not yield anything; it is possible that adjoints of degree $d-4$ provide tighter results, but in that case the conditions imposed by the definition of adjoint are not linearly independent in general (which is good; if they were, the dimension of the adjoint space would be $\gg=d+1$ which precludes the solution of any interesting case). The condition for having 4 or less intersection points in $\pj^3(\K)\setminus \pj^3(\F)$ translates into $\gg\leq6$ which is very suggestive.

\section*{Acknowledgements}

The authors thank Gian Pietro Pirola for his useful remarks, specially for pointing out the work of O. Zariski on this topic.

\bibliography{biblio}
\bibliographystyle{plain}

\end{document}